\documentclass[12pt]{article}
\usepackage{fullpage,indentfirst}
\usepackage{amsfonts,amssymb,stmaryrd}
\usepackage{diagrams}

\newarrow{mapsto}|--->

\newtheorem{theo}{Theorem}

\newtheorem{lemma}{Lemma}[section]

\newcounter{defini}[section]
\newcommand\defi{
\refstepcounter{defini}
\noindent{\bf Definition\  \thesection
    .\arabic{defini}\ }
}

\def\RR{\mathbb{R}}
\def\CC{\mathbb{C}}
\def\ZZ{\mathbb{Z}}
\def\NN{\mathbb{N}}

\def\longto{\longrightarrow}

\def\a{\alpha}
\def\om{\omega}
\def\proof{\noindent{\bf Proof~:~}}
\def\carre{\hfill $\square$}

\def\vsp{\vspace{3mm}}
\def\proof{\noindent{\bf Proof~:~}}
\def\carre{\hfill $\square$}

\newcounter{paragraf}[section]
\renewcommand{\theparagraph}{\thesection.\arabic{paragraf}}
\renewcommand{\paragraph}
{\refstepcounter{paragraf}
\vspace{3mm}
{\bf \theparagraph}\hspace{0.2em}--- 
}

\begin{document}

~\vspace{3cm}
\begin{center}
  {\bf \Large Balanced configurations of $2n+1$ plane vectors}\\

\vsp
N. Ressayre

\end{center}

\vspace{1cm}
\section{Introduction}

A plane configuration $\{v_1,v_2,\cdots,v_m\}$ (where $m$ is a
positive  integer) of vectors of $\RR^2$ is said to be {\it balanced}
if for any index  
$i\in \{1,\cdots,m\}$ the multiset
$$
\{\det(v_i,v_j)\;:\;j\neq i\}
$$
is symmetric around the origin. 
A plane configuration is said to be {\it uniform} if
every pair of vectors is linearly independent. 

E. Cattani, A. Dickenstein and B. Sturmfels introduced this notion in
\cite{CatDickSt,CatDick} for its relationship with multivariable
hypergeometric functions in the sense of Gel'fand, Kapranov and
Zelevinsky (see \cite{GKZhypergeom,GKZeulerintegral}). 

Balanced plane configurations with at most six vectors have been
classified in \cite{CatDickSt}. 
With the help of computer calculation, E. Cattani, A. Dickenstein
classified the balanced plane configurations of seven vectors in
\cite{CatDickSt}. 
Moreover, they conjectured that any uniform balanced plane
configuration is ${\rm GL}_2(\RR)$-equivalent to a regular $(2n+1)$-gon
(where $n$ is  a positive integer). 
In this note, we prove this conjecture.

\section{Statement of the result}

Let $m$ be a positive integer. 
 
\defi
A configuration $\{ v_1,\cdots,v_m\}$ is said to be ${\it balanced}$
if for all $i=1,\cdots,m$ and for all $x$ in $\RR$ the cardinality of
the set $\{j\neq i\;:\;\det(v_i,v_j)= x\}$  equals those of
the set $\{j\neq i\;:\;\det(v_i,v_j)= -x\}$.

\vsp
\defi 
A  balanced configuration $\{v_1,\cdots,v_m\}$ is said
to be {\it uniform} if for any pair $i\neq j$, the  vectors $v_i$,
$v_j$ are linearly independent. 

\vsp
\noindent
{\bf Remark:}
Assume $\{v_1,\cdots,v_m\}$ is balanced and $m$ even. 
Then, the multiset $\{\det(v_1,v_j)\;:\;j=2,\cdots,m\}$ is symmetric
around 0 and of odd cardinality; so it contains $0$. 
Then, $\{v_1,\cdots,v_m\}$ is not uniform.
From, now on we are only interested in configurations with an odd
number of vectors. So, we assume that $m=2n+1$ for an integer $n$.

\vsp 
Let us identify $\RR^2$ with the field $\CC$ of complex numbers. 
To avoid any confusion with index-numbers, we denote by 
$\sqrt{-1}$ the complex number $i$.
Denote by ${\mathbb U}_m$ the set of $m^{\rm th}$-roots of 1. 

Set $\om=e^{\frac{2\sqrt{-1}\pi}{m}}$. 
Then, ${\mathbb U}_m=\{w^k\;:\;k=0,\cdots 2n\}$. 
For all integers $k$ and $a$, we have 
\begin{eqnarray}
  \label{eq:kpml}
  \det(\om^k,\om^{k+a})=-\det(\om^k,\om^{k-a}).
\end{eqnarray}

In particular,  ${\mathbb U}_m$ is uniform balanced configuration.

One can note that the group ${\rm GL}_2(\RR)$ acts naturally on the set 
of  balanced (resp. uniform balanced)  configurations of $m$ vectors. 
Indeed, if $g\in {\rm GL}_2(\RR)$  then 
$\det(g.v_i,g.v_j)=\det(g)\det(v_i,v_j)$.

The aim of this note is to prove the 

\begin{theo}
 \label{th}
For any odd integer $m$,  ${\rm GL}_2(\RR)$  acts transitively  on the set 
of  uniform balanced  configurations of $m$ vectors. 
\end{theo}

In other words, modulo ${\rm GL}_2(\RR)$, ${\mathbb U}_m$ is the only 
uniform balanced configuration of $m$ vectors. 

\section{The proof}

\vsp
\paragraph
Let us fix some notation and convention. 
The set $\{0,\cdots,2n\}$ is denoted by $I$. 

\vsp
\defi
Let us recall that we identify $\RR^2$ with the field $\CC$ of complex numbers. 
Let $\{v_0,\cdots,v_{m-1}\}$ be a uniform configuration of $m$ points
in $\RR^2$. Each $v_i$ has a unique polar form $v_i=\rho_i e^{\a_i}$
with $\rho_i$ in $]0;+\infty[$ and $\a_i$ in $[0;2\pi[$. 
The set $\{v_0,\cdots,v_{m-1}\}$ is said to be {\it labeled by increasing arguments} if 
$$
\a_0<\a_1<\cdots<\a_{m-1}.
$$

\vsp
\noindent
{\bf Convention 1}
Let $i\in I$. For all $k$ in $\ZZ$ which equals
$i$ modulo $m$, we also denote by $v_k$ the vector $v_i$. 

\vsp
The first step of the proof is to show that any uniform configuration satisfy 
equations similar to Equations (\ref{eq:kpml}).  Precisely, we have:

\begin{lemma}
\label{lem:ordre}
Let  ${\cal C}=\{v_0,\cdots,v_{2n}\}$ be a uniform balanced configuration 
labeled by increasing arguments.
Then,\\
$$
\det(v_k,v_{k+a})=-\det(v_k,v_{k-a})\hspace{2cm} \forall k,\,a\in\ZZ
$$
\end{lemma}

\proof
We denote by  ${\cal  P}_2(I)$ the set of pairs of elements of $I$. 
The fact that  ${\cal C}$  is uniform balanced can be
formulated as follow. 
For all $i\in I$, there exists a part ${\cal  P}_2^i(I)$ of 
${\cal  P}_2(I)$ such that:
\begin{itemize}
\item $I-\{i\}$ is the disjoint union of the elements of ${\cal  P}_2^i(I)$, and 
\item $\forall \{k,l\}\in{\cal  P}_2^i(I)\hspace{1cm} 
\det(v_i,v_k)=-\det(v_i,v_l)\neq 0 $.
\end{itemize}

For any pair $\{k,l\}\in{\cal  P}_2(I)$, the set of
vectors $v\in\RR^2$ such that $\det(v,v_k)=-\det(v,v_l)$ is the
vectorial line generated by $v_k+v_l$ (let us recall that $v_k,v_l$ are linearly independent). 
In particular, since ${\cal C}$ is uniform there
exists at most one $i\in I$ such that $\det(v_i,v_k)=-\det(v_i,v_l)$.
This means that for any $i\neq j$ the set 
${\cal  P}_2^i(I)\cap {\cal  P}_2^{j}(I)$ is empty.

Moreover, the cardinality of ${\cal  P}_2^i(I)$ equals $n$ for all
$i\in I$. Then, the cardinality of $\bigcup_{i\in I}{\cal  P}_2^i(I)$
equals  $nm$, that is the cardinality of ${\cal  P}_2(I)$. It follows that
$
\bigcup_{i\in I}{\cal  P}_2^i(I)={\cal  P}_2(I)
$.
In other words, there exists a map
$$
\begin{array}{r@{}ccc}
\phi\;:\;&{\cal  P}_2(I)&\longto&I,
\end{array}
$$ 
such that, for all $\{k,l\}\in {\cal  P}_2(I)$, we have :
$$
\det(v_{\phi(\{k,l\})},v_k)=-\det(v_{\phi(\{k,l\})},v_l).
$$ 

It is sufficient to prove the lemma for
$a=-n,\cdots,-1,1,\cdots,n$; and by symmetry for 
$a=1,\cdots,n$. 
We prove this by decreasing induction going from $a=n$ to $a=1$. 

Assume $a=n$ and fix $k$. Relabeling the vectors, we may assume that
$k=n+1$. Then, we have to prove that :
$\det(v_{n+1},v_0)=-\det(v_{n+1},v_{1})$, that is, $\phi(\{0,1\})=n+1$.

Note that the set of $i\in I$ such that $\det(v_0,v_i)$ is positive
(that is, such that $\a_i-\a_0<\pi$) is of cardinality $n$. 
Then, by Convention 1 $\a_n-\a_0<\pi$. 

For all $t=0,\cdots,n-1$,  since $v_{\phi(\{t,t+1\})}$ belongs to
$\RR(v_t+v_{t+1})$, its argument $\a_{\phi(\{t,t+1\})}$ belongs to
$]\pi+\a_t;\pi+\a_{t+1}[$. In particular, each one of the $n$ intervals   
$]\pi+\a_t;\pi+\a_{t+1}[$ (for $t=0,\cdots,n-1$) contains one of the
$\a_i$ for $i=n+1,\cdots,2n$. So, $\a_{\phi(\{0,1\})}$ is the only
$\a_i$ in the interval $]\pi+\a_0;\pi+\a_{1}[$. It follows that
$\phi(\{0,1\})=n+1$.  

Suppose now the proposition proved for $a=n,\cdots,n-u+2$ 
(with $n\geq u\geq 2$) and prove that it is true for $a=n-u+1$. 
As before, it is sufficient to prove that :
$$
\begin{array}{l@{}ll}
\phi(\{0,u\})&=\frac{u}{2} &{\rm if\ } $u$ {\rm \ is\ even}\\[1mm]
             &=\frac{u+m}{2}=\frac{u+1}{2}+n&  {\rm if\ } $u$ {\rm \ is\ odd}
\end{array}
$$

Since, $v_{\phi(\{0,u\})}$ belongs to $\RR.(v_0+v_u)$, we have : 
$$
\phi(\{0,u\})\in\{1,\cdots,u-1\} \cup \{n+1,\cdots,n+u\}. 
$$

Let us assume that $u=2v$ is even. 
For $w=0,1,\cdots v-1$, 
we have $\phi(\{0,2w+1\})=n+1+w$. But, two  elements of 
${\cal P}_2^{n+1+w}$ are disjoint. 
So, $\phi(\{0,u\})\not\in \{n+1,\cdots,n+v\}$.
In the same way, for $w=1,\cdots,v-1$, we have : $\phi(\{0,2w\})=w$. 
And so,  $\phi(\{0,u\})\not\in \{1,\cdots,v-1\}$.
For $w=0,1,\cdots v-1$, we have $\phi(\{u,u-2w-1\})=n+u-w$. Then, 
$\phi(\{0,u\})\not\in \{n+v+1,\cdots,n+u\}$.
For $w=1,\cdots,v-1$, we have $\phi(\{u,u-2w\})=u-w$. Then, 
 $\phi(\{0,u\})\not\in \{v+1,\cdots,u-1\}$.

Finally, the only possible value for  $\phi(\{0,u\})$ is $v$.

The proof is analog if $u=2v+1$ is odd.
\carre

\vsp
Lemma~\ref{lem:ordre} has a very useful consequence:

\begin{lemma}
\label{lem:aetA}
We keep notation of Lemma~\ref{lem:ordre}. 
We also use Convention 1.

Then, for all $k=0,\cdots,2n$ we have :
$$
\det(v_k,v_{k+1})=\det(v_0,v_1),
$$
and
$$
\det(v_k,v_{k+n})=\det(v_0,v_n).
$$
\end{lemma}

\proof
Lemma~\ref{lem:ordre} shows that for all integer $k$ we have 
$
\det(v_k,v_{k+1})=\det(v_{k+1},v_{k+2}).
$
The first assertion follows immediately.

For all $k$, we also have 
$
\det(v_k,v_{k+n})=\det(v_{k+n},v_{k+2n}).
$
Since $n$ is prime with $m=2n+1$, this implies the second assertion.
\carre

\paragraph
Let ${\cal C}=\{v_0,\cdots,v_{2m}\}$ be a uniform balanced
configuration labelled by increasing arguments. We are going to prove

\vsp
\noindent
{\bf Claim 1:} $v_0,v_n$ and $v_{n+1}$ determine ${\cal C}$.

\vsp
Indeed, we are going to construct successively
$v_1,v_{n+2},v_2,v_{n+3},v_3,v_{n+4}\ldots$. Set
$A_1:=\det(v_n,v_{n+1})$ and $A_n:=\det(v_0,v_n)$. 
Assume that we have constructed $v_1,v_{n+2},\cdots,v_{i-1},v_{n+i}$
(for $1\leq i\leq n-1$). By Lemma~\ref{lem:aetA}, we have:\\

$~\hfill\hfill
\det(v_{i-1},v_i)=A_1\hfill{\rm and}\hfill
\det(v_i,v_{n+i})=A_n. 
\hfill(2)\hfill
$\\

Then,
$$
v_i=\frac{A_1}{\det(v_{i-1},v_{n+i})}v_{n+i}+
\frac{A_n}{\det(v_{i-1},v_{n+i})}v_{i-1}.
$$
But, since by Convention 1, $v_{n+i+n}=v_{i-1}$, we have:
$\det(v_{i-1},v_{n+i})=-A_n$. 
Finally, we obtain:
$$
v_i=\frac{A_1}{A_n}v_{n+i}-v_{i-1}.
$$

In the same way, using\\ 

$~\hfill\hfill
\det(v_{n+i},v_{n+i+1})=A_1\hfill{\rm and}\hfill
\det(v_{n+i+1},v_i)=A_n; 
\hfill(3)\hfill
$

\noindent
we obtain:
$$
v_{n+i+1}=-\frac{A_1}{A_n}v_{i}-v_{n+i}.
$$
Claim 1 follows.

\paragraph{}
Inspired by the proof of Claim 1, we define two sequences of vectors
of $\RR^2$ (with a parameter $t\in \RR$) as follows.

Start with \\

$~\hfill\hfill
U=
\left (
  \begin{array}{c}
1\\0
  \end{array}
\right )
\hfill
V=\left (
  \begin{array}{c}
0\\1
  \end{array}
\right )
\hfill
w_0(t)=\left (
  \begin{array}{c}
t\\-1
  \end{array}
\right ).
\hfill\hfill
$
  
Set $A=\det(V,w_0)=-t$ and note that $\det(U,V)=1$. 
Then we define $w_i(t)$ and $u_i(t)$ by induction:
$$\left\{
\begin{array}{ll}
w_0(t)&{\rm is\ already\ defined}\\
u_0(t)=U\\
u_{i+1}(t)=-tw_i(t)-u_i(t)\\
w_{i+1}(t)=tu_i(t)-w_i(t)
\end{array}
\right.
$$

\paragraph{}
\label{par:deftC}
Let ${\cal C}=\{v_0,\cdots,v_{2m}\}$ be a uniform balanced
configuration labelled by increasing arguments. 
Then, there exits a unique $g_{\cal C}\in {\rm GL}_2(\RR)$ such that $g_{\cal C}.v_0=U$
and $g_{\cal C}.v_n=V$. Since 
$\det(v_0,v_n)=-\det(v_0,v_{n+1})$  (see Lemma~\ref{lem:aetA}), 
there exists a unique $t_{\cal C}\in\RR$ such that
$g_{\cal C}.v_{n+1}=w_0(t_{\cal C})$. 
Then, the proof of Claim 1 implies

\begin{lemma}
\label{lem:suiteC}
With above notation, for all $i=0,\cdots,n-1$, we have:
$$
g_{\cal C}.v_{n+i+1}=w_i(t_{\cal C})\hspace{2cm}{\rm and}\hspace{2cm}
g_{\cal C}.v_i=u_i(t_{\cal C}).
$$  
Moreover, $w_n(t_{\cal C})=U$ and $v_n(t_{\cal C})=V$.
\end{lemma}

\paragraph
Now, we are interested in the equation $w_n(t)=U$. 

Useful properties of the functions $t\mapsto u_i(t)$ and 
$t\mapsto w_i(t)$ are stated in

\begin{lemma}
\label{lem:polpair}
Denote by $(x^*,y^*)$ the coordinate forms of $\RR^2$. 
Then, for all $i\geq 1$, we have:
\begin{enumerate}
\item $x^*(u_i(t))$ is an even  polynomial function of degree $2i$,
\item $y^*(u_i(t))$ is an odd polynomial function of degree $2i-1$,
\item $x^*(w_i(t))$ is an odd polynomial function of degree $2i+1$, and
\item $y^*(w_i(t))$ is an even polynomial function of degree $2i$.
\end{enumerate} 
In particular, the equation $w_n(t)=U$ has at most $n$ solutions.  
\end{lemma}

\proof
The proof of the four assumptions is an immediate induction on $i$.

We can note that $y^*(w_n(0))\neq 0$. 
Then, by Assertion $(iv)$, the equation $y^*(w_n(t))=0$ has at most $2n$ solutions :
$-t_j<\cdots<-t_1<t_1<\cdots<t_j$ (with $j\leq n$).
Since $x^*(w_n(t))$ is an odd polynomial function, at most one element
of a pair $\pm t_k$ is a solution of the equation $x^*(w_n(t))=1$.
This ends the proof of the lemma.
\carre

\paragraph
\label{par:deftk}
Our goal is now to construct geometrically $n$ solutions of the
equation $w_n(t)=U$.

Let me recall that we have identified $\RR^2$ with $\CC$. 
Consider 
${\mathbb U}_m=\{\om^i\,:\,i=0,\cdots 2n\}$.  
Let us fix $k\in\{1,\cdots n\}$.

Denote by $g_k$ the element of ${\rm GL_2}(\RR)$ such that $g_k.1=U$ and 
$g_k.\om^k=V$. 
Let $t_k$ be the unique real number such that
$g_k.\om^{-k}=w_0(t_k)$.
Explicitly, $t_k=\frac{1}{\sin(2k\pi/m)}$. 

For all $i\in\ZZ$, we have :\\

$~\hfill
\det(\om^{-2k(i-1)},\om^{-2ki})=\det(\om^k,\om^{-k})
\hfill
\det(\om^{-2ki},\om^{-2k(n+i)})=\det(\om^0,\om^k),
\hfill
$\\

and\\

$~\hfill
\det(\om^{-2k(n+i)},\om^{-2k(n+i+1)})=\det(\om^k,\om^{-k})
\hfill
\det(\om^{-2k(n+i+1)},\om^{-2ki})=\det(\om^0,\om^k).
\hfill
$

\vsp
Then, the sequence $(g_k.\om^{-2ki})_{i\in\NN}$ satisfies 
Relations $(2)$ and $(3)$, with 
$A_1=\det(V,w_0(t_k))$ and $A_n=\det(U,V)$. 
This implies that 

$$
\hspace{4.5cm}
w_i(t_k)=g_k.\om^{-k(1+2i)}, 
{\rm \ for\  all\ }i\geq 0\hspace{4.5cm}(4).
$$

In particular, $t_k$ satisfies $w_n(t_k)=U$. 

With Lemma~\ref{lem:polpair}, this implies the

\vsp
\begin{lemma}
\label{lem:eq}
We have:
$$
\{t\in\RR\;:\;w_n(t)=U\}=\{\frac{1}{\sin(2k\pi/m)}\;:\;k=1,\cdots n\}.
$$  
\end{lemma}

\paragraph
\noindent
{\bf Proof of Theorem \ref{th}}
Let ${\cal C}=\{v_0,\cdots,v_{2n}\}$ be a uniform balanced configuration 
labeled by  increasing arguments.  
We define $g_{\cal C}\in{\rm GL}_2(\RR)$ and $t_{\cal C}\in\RR$ as 
in Paragraph~3.\ref{par:deftC}.
Then, by Lemmas~\ref{lem:suiteC} and \ref{lem:eq}, there exists a unique 
$k_{\cal C}=1,\cdots n$ such that $t_{\cal C}=\frac{1}{\sin(2k_{\cal C}\pi/m)}$. 
Let $g_{k_{\cal C}}\in{\rm GL}_2(\RR)$ defined as in Paragraph~3.\ref{par:deftk}.

Then, by Lemma~\ref{lem:suiteC} and Equalities~$(4)$, we have:

\begin{diagram}
v_0&\rmapsto^{g_{\cal C}}&U&\rmapsto^{g_{k_{\cal C}}^{-1}}&1\\
 v_n&\rmapsto^{g_{\cal C}}&V&\rmapsto^{g_{k_{\cal C}}^{-1}}&\om^{k_{\cal C}}\\
 v_{n+i+1}&\rmapsto^{g_{\cal C}}&w_i(t_{\cal C})&
\rmapsto^{g_{k_{\cal C}}^{-1}}&\om^{-k(1+2i)}&
{\rm\hspace{1cm} \ for\ all\ }i=0,\cdots,n-1\\
 v_i&\rmapsto^{g_{\cal C}}&v_i(t_{\cal C})&\rmapsto^{g_{k_{\cal C}}^{-1}}&\om^{-2ki}&
{\rm \hspace{1cm} \ for\ all\ }i=0,\cdots,n-1
\end{diagram}

Theorem~\ref{th} follows.
\carre

\bibliographystyle{smfalpha}
\bibliography{biblio}

\begin{center}
  -\hspace{1em}$\diamondsuit$\hspace{1em}-
\end{center}

\vspace{5mm}
\begin{flushleft}
Nicolas Ressayre\\
Universit\'e Montpellier II\\
D\'epartement de Math\'ematiques\\
Case courrier 051-Place Eug\`ene Bataillon\\
34095 Montpellier Cedex 5\\
France\\
e-mail:~{\tt ressayre@math.univ-montp2.fr}  
\end{flushleft}

\end{document}